\documentclass[13pt]{amsart}
\usepackage{graphicx}
\usepackage{enumerate}
\newtheorem{thm}{Theorem}[section]

\newtheorem{lem}[thm]{Lemma}

\theoremstyle{definition}

\theoremstyle{remark}

\numberwithin{equation}{section}

\def\d{{\rm d}}

\def\tilda{\widetilde}
\def\midd{\,|\,}
\def\E{{\rm E}}
\def\arr{\rightarrow}
\def\sumin{\sum_{i=1}^n}
\def\hop{\smallskip\noindent}

\begin{document}
\parindent20pt
\baselineskip17pt

\def\today{May 2005}

\title{Bayesian Bivariate Survival Estimation}
\author{J.K.~Ghosh, N.L.~Hjort, C.~Messan and R.V.~Ramamoorthi}

\begin{abstract}
There is no easy extension of Kaplan--Meier 
and Nelson--Aalen to the bivariate case, 
and estimating bivariate survival distributions 
nonparametrically is associated with various 
non-trivial problems. The Dabrowska estimator
will for example associate negative mass to some
subsets. Bayesian methods hold some promise as 
they will avoid the negative mass problem, but
are also prone to difficulties. We simplify and extend
an example by Pruitt to show that posterior distribution
from a Dirichlet process prior is inconsistent.
We construct a different nonparametric prior 
via Beta processes and provide an updating scheme
that utilizes only the most relevant parts of the likelihood, 
and show that this leads to a consistent estimator. 


\end{abstract}

\date{\today}

\maketitle
\section{Introduction}

\noindent
Bivariate survival or waiting times to occurrence of some event
has several interesting applications. For example a husband and
his wife may both be exposed to a common risk and one might want
to know if males and females have the same survival distribution.
Another situation might be the study of effects of active and
passive smoking when only one of the spouses smokes.

In the univariate case the most popular nonparametric
estimator is the Kaplan--Meier one. The bivariate case
remains surprisingly difficult in spite of a lot of work.
Dabrowska (1988)
constructs a bivariate analog of the Kaplan--Meier estimate which
is consistent but is not a proper survival distribution in that it
assigns negative mass to some events. The same is true of an
earlier estimate of Langberg and Shaked (1982). 
Dabrowska (1988)
mentions another estimate due to Bickel which is consistent and avoids
negative masses but does not seem to utilize all of the data.
Pruitt (1991)
shows mathematically that the negative mass
assigned by Dabrowska's estimator can be substantial even though
the estimate assigns positive mass to all upper orthants.
A fairly complete survey of some of these approaches
can be found in Andersen, Borgan, Gill and Keiding (1993, Ch.~X)

Pruitt (1988)
also shows that the Bayes estimate of the
bivariate survival function with a Dirichlet process prior can be
dramatically inconsistent. In Pruitt (1991) he suggests an
independence assumption weaker than Dabrowska's and explains
how in this context one can construct a new bivariate
Kaplan--Meier estimate which is a proper survival function.

The univariate Kaplan--Meier estimate can be obtained as a limit of
Bayes estimates with respect to Dirichlet process priors,
more generally Beta processes (see Hjort, 1990)
and even more generally through neutral to the right priors
(Hjort, op.~cit.).
In this paper we develop a natural generalization of Beta processes
for the bivariate case and derive an estimator using an essentially
Bayesian approach. The qualification `essentially' is used since
we do not use the full likelihood.

The organization of this paper is as follows. In Section 2 we
discuss Pruitt's example and give a simple proof of inconsistency.
In Section 3 we consider discrete distributions for the survival
time and through a convenient parametrization relate it to the
distribution of the censored data. The discussion in this section
provides some insight into the difficulties associated with the
bivariate case.
Then in section~4 we look at the time-discrete case and develop the
Bayesian analysis. The last section considers the general case.
Though we restrict ourselves to Beta processes 
much of our arguments and results would go through
for neutral to the right priors, but we do not pursue this here.

\section{Pruitt's example}

\noindent
We follow Pruitt (1988) but provide a rather easy proof of
inconsistency. Let $(T_1,T_2)$ be two survival times,
$(C_1,C_2)$ the censoring variables,
and let $Z_j=T_j\wedge C_j$ and $\Delta_j=I\{T_j\le C_j\}$
for $j=1,2$ be the observed random variables.

Let $(C_1,C_2)$ be discrete with mass 1/3 at each of
$(1,3),(3,1),(4,4)$. The $(T_1,T_2)$ pair is assumed
to take values in $\mathcal{T}= [1,3]\times [1,3]$.
The prior is a Dirichlet process with parameter $M\alpha$
where $\alpha$ is the uniform distribution
on $[1,3]\times [1,3]$. We will show that if the true
distribution $P_0$ is uniform on
$A=[1,2]\times[1,2]\cup [2,3]\times[2,3]$
then the posterior is inconsistent. Much more delicate
calculations are needed to show that the estimate converges weakly
to $(1/3)P_{0}+(2/3)\alpha$.
Let $(\tilda{Z},\tilda{\Delta})$
denote the observations
$(Z_{1,i},\Delta_{1,i}),(Z_{2,i},\Delta_{2,i})$
for $i=1,\ldots,n$, with
$$\tilda Z_i=(Z_{1,i},Z_{2,i})
  \quad {\rm and} \quad
  \tilda\Delta_i=(\Delta_{1,i},\Delta_{2,i}). $$

Consider $B=[1,2]\times[2,3]$. Then $P_{0}(B)=0$. Let
$\E\{P(B)\midd(\tilda{Z},\tilda{\Delta})\}$ be the
Bayes estimator. We shall show that this does not go 0 in
probability as $n\arr \infty$.

Recall that if $E$ is a subset of $\mathcal{T}=[1,3]$ then with a
Dirichlet($\alpha$) prior for the distributions on $\mathcal{T}$,
the posterior given $E$ is the mixture
\[\int D_{\alpha+\delta(\tilda t)}\,\alpha(\d \tilda{t}\midd E), \]
where $\delta(\tilda t)=\delta_{\tilda t}$
denotes unit mass at point $\tilda t$,
i.e.~$\delta(\tilda t)(B)=I\{\tilda t\in B\}$,
and $\alpha(\cdot\midd E)$ is the conditional distribution of
$\tilda{t}$ given $E$. Extending this to $n$  observations,
and thinking of the censored observations as specifying subsets of
$\mathcal{T}$, for instance $(\Delta_1,\Delta_2)=(0,1)$
as corresponding to the horizontal line at $Z_{1}$,
the posterior is a mixture of the form
\[\int D_{\alpha+\sum_{i=1}^n\delta(\tilda t_i)}\,
  \alpha(\d\tilda t_i\midd(\tilda Z_1,\tilda\Delta_1),
   \ldots,(\tilda Z_n,\tilda\Delta_n)). \]
Consequently the Bayes estimate of $P$ is
\[\frac{M}{M+n}\alpha
+\frac{n}{M+n}\frac{1}{n}\sumin
   \E\{\delta(\tilda t_i)\midd
  (\tilda{Z},\tilda{\Delta})\}. \]
Since $P_0$ is continuous we may assume that
$\tilda t_1,\ldots,\tilda t_n$ are all distinct.
Given this assumption, $\tilda t_1,\ldots,\tilda t_n$
are independent with common distribution $\alpha$.
Hence the Bayes estimate of $P(B)$ is
\[\frac{M}{M+n}\alpha(B)
+\frac{n}{M+n}\frac{1}{n} \sumin
  \E\left[\delta(\tilda t_i)(B)\midd((Z_{1,i},\Delta_{1,i}),
  (Z_{2,i},\Delta_{2,i}))\right]. \]
The expectation in the above expression is calculated taking
$\tilda t_1,\ldots,\tilda t_n$ to be independently
drawn from $\alpha$. As $n\arr\infty$, the first term in
the last expression goes to 0 so we shall concentrate on the
second term.

Note that under the model it is not possible to have both
$\Delta_1$ and $\Delta_2$ equal to 0. Further, if
$\Delta_{1,i}=1,\Delta_{2,i}=1$ for some $i$
then the conditional distribution of $\tilda t_i$
given $(Z_{1,i},\Delta_{1,i}),(Z_{2,i},\Delta_{2,i})$
is degenerate at $(Z_{1,i},Z_{2,i})$. So observations with
$\Delta_{1,i}=1,\Delta_{2,i}=1$ contribute to the Bayes estimate
of $P(B)$ only if the corresponding $(Z_{1,i},Z_{2,i})$ is in $B$.
Since $P_0(B)=0$, under $P_0$ there will be no such
observations and so will not contribute to the Bayes estimate
of $P(B)$. Thus the only case that we need to consider is
when exactly one of $\Delta_{1,i}$ and $\Delta_{2,i}$ is 1
and the other is 0. Thus
\vskip-20pt
\begin{setlength}{\multlinegap}{0pt}
\begin{multline*}
\sumin\E\{\delta_{\tilda t_i}(B)\midd(Z_{1,i},\Delta_{1,i}),
   (Z_{2,i},\Delta_{2,i})\} \\
=\sum_{i:\Delta_{1,i}=0,\Delta_{2,i}=1}
   \E\{\delta_{\tilda{t_{i}}}(B)|(Z_{1,i},\Delta_{1,i}\},
   (Z_{2,i},\Delta_{2,i})) \\
+\sum_{i:\Delta_{1,i}=1,\Delta_{2,i}=0}
   \E\{\delta_{\tilda{t_{i}}}(B)|(Z_{1,i},\Delta_{1,i}),
   (Z_{2,i},\Delta_{2,i})\}.
\end{multline*}
\end{setlength}


Next consider an $i$ for which
$Z_{1,i}=z_{1,i},Z_{2}=z_{2,i}$
and $\Delta_{1,i}=0,\Delta_{2,i}=1$.
In this case, we know that $T_{1} \in [1,3]$ and $ T_{2}=z_{2,i}$.
Consequently, under $\alpha$, the conditional distribution
of $T_1$ given
$\{Z_{1,i}=z_{1,i},\Delta_{1,i}=0,Z_{2}=z_{2,i},\Delta_{2,i}=1\}$
is uniform on $[1,3]$. From these observations it follows that
if $Z_{2,i}\in [2,3]$,
\[\E\{\delta_{\theta_i}(B)\midd Z_{1,i},\Delta_{1,i}=0,Z_{2,i},
  \Delta_{2,i}=1\}
= \E\{I_{B}(T_1,T_2)\midd T_{1,i}\geq 0,T_{2,i}=Z_{2,i}\} \]
and equal to zero otherwise.
Now if $Z_{2,i}\in [2,3]$ then
\[\E\{I_B(T_1,T_2)\midd T_{1,i}\geq 0,T_{2,i}=Z_{2,i}\}
  =\frac{1}{2}\int_{1}^{2}\d t_1=\frac{1}{2}. \]
One handles the case of $\Delta_{1,i}=1,\Delta_{2,i}=0$ similarly.
Hence
\begin{multline*}
\frac{n}{M+n}
\frac{1}{n}\sumin \E\{\delta_{\tilda{t}_{i}}\midd
  (Z_{1,i},\Delta_{1,i}),(Z_{2,i},\Delta_{2,i})\} \\
= \frac{n}{M+n}\Bigl\{\frac{1}{2n}
   \sum_{i\colon \Delta_{1,i}=0,\Delta_{2,i}=1}I_{[2,3]}(Z_{2,i})
  +\frac{1}{2} \sum_{i\colon
  \Delta_{1,i}=1,\Delta_{2,i}=0}I_{[1,2]}(Z_{1,i})\Bigr\}
\end{multline*}
and as $n\arr \infty$ this goes to
\[
(1/2)\bigl[ P\{Z_2\in [2,3],\Delta_1=0,\Delta_2=1\}
  +P\{Z_1\in[1,2],\Delta_1=1,\Delta_2=0\} \bigr]. 
\]
Using the fact that
$$P\{Z_2\in[2,3],\Delta_1=0,\Delta_2=1\}
  =P\{Z_1\in[1,2],\Delta_1=1,\Delta_2=0\}=(1/2)(1/3), $$
the limit is seen to be 1/6. This shows that the Bayes estimator
is inconsistent and hence so is the posterior. \qed

\smallskip
It might appear that the difficulty stems from $P_0$ and the prior
guess not having the same support. However, a little reflection
will show that the inconsistency phenomenon will occur 
even if $P_0$ and the prior guess have the same support.

\section{Parametrization of distributions}

\noindent
Before we proceed to the details it is convenient to have an
overview of the developments in this section.

Consider the usual one-dimensional case, where $T$ is a life-time
with unknown distribution $F$ and  $C$, independent of $T$, is a
censoring variable with a fixed known distribution $G$. The
observation is $(Z,\Delta)$ where $Z = T\wedge C$ and
$\Delta=I\{T\leq Z\}$. Let $\mathcal{M}(T), \mathcal{M}(Z,\Delta)$
denote respectively the set of all distributions of $T$ and $(Z,\Delta)$.
Under minimal assumptions the map
$\phi\colon\mathcal{M}(T) \mapsto \mathcal{M}(Z,\Delta)$,
which takes any $F$ to $F^{(Z,\Delta)}$,
is one-one and onto. In particular, Peterson (1977)
has shown that the inverse of the empirical
distribution leads to the Kaplan--Meier estimate.

The situation is drastically different in the two- and higher-dimensional
case. Here if
$\mathcal{M}(\tilda{T}),\mathcal{M}(\tilda{Z},\tilda{\Delta})$
denote respectively the set of all distributions of $\tilda{T},
(\tilda{Z},\tilda{\Delta})$, the corresponding map
$\phi\colon\mathcal{M}(\widehat{T})\mapsto
\mathcal{M}(\tilda{Z},\tilda{\Delta})$, which takes any
$F$ to $F^{(\tilda{Z},\tilda{\Delta})}$,
is one-one but far from onto. In fact many empirical
distribution functions fall outside the range.
What we do in this section is, in an abstract sense,
to define an `inverse' to all of
$\mathcal{M}(\tilda{Z},\tilda{\Delta})$, and we do this
through a reparametrization of the underlying distributions.

We first reparametrize the distribution of $T$,
then carry out a similar operation for the distribution of
$(\tilda{Z},\tilda{\Delta})$ and show that using only part
of the distribution of $(\tilda{Z},\tilda{\Delta})$  we
can recover the distribution of $T$. This makes it transparent
that the function $\phi$ is not onto.
We then use the reparametrization to write the likelihood
function. We note that it contains a component that is complex
and propose that it be ignored. It will be seen in the next
section that working with this incomplete likelihood
enables one to develop a tractable posterior.

Let $\tilda{T}=(T_{1},T_{2})$ have the distribution $P$,
a discrete probability measure on the positive
quadrant, for concreteness and without essential loss
of generality taken to have support in
$\{1,2,3,\ldots\}\times\{1,2,3,\ldots\}$.
We begin by parametrizing $P$ as follows.
Let first $T^{*}=T_{1}\wedge T_{2}$, and let
$\epsilon=0$ if $T_1=T_2$;
$\epsilon=1$ if $T_1>T_2$;
and $\epsilon=2$ if $T_1<T_2$.
Furthermore, let
$$
\begin{cases}
    P^{*} & \text{distribution of } T^{*}, \\
     P^{\epsilon}(\cdot\midd T^{*}) & \text{distribution of }
   \epsilon \text{ given } T^{*}, \\
P_{1}(\cdot\midd T^{*},1) & \text{distribution of }T_1
  \text{ given } T^{*},\epsilon=1, \\
P_{2}(\cdot\midd T^{*},1) & \text{distribution of }T_{2}
  \text{ given } T^{*},\epsilon=2.
  \end{cases} 
$$
To incorporate the censoring model, let $\tilda{C}=
(C_{1},C_{2})$ be a pair of censoring variables, independent of
$\tilda{T}$ and with distribution $G$. We assume that for
every $(x,y)$ in $\mathbb{R}\times \mathbb{R}$ there is an
$(x', y')$ with positive $G$-probability such that
$x'>x,\,y^\prime >y$.

As before we only observe
\[( \tilda{Z},\tilda{\Delta})
   =(Z_1,\Delta_1,Z_2,\Delta_2), \]
where $Z_j=T_j\wedge C_j$ and $\Delta_j=I\{T_j\le C_j\}$.
We  parametrize the distribution $P^{Z}$ of
$(\tilda Z,\tilda\Delta)$ as follows:
$$
\begin{cases}
Z^{*}= Z_1\wedge Z_2 = (T_1\wedge T_2) \wedge
(C_1\wedge C_2), \\
\Delta^{*}= I\{T^{*}\leq C^{*}\} &\text{where } C^{*}=C_{1}\wedge
C_2 \\
\eta = 0  & \text{if } Z_{1}=Z_{2}, \\
\eta =1  & \text{if } Z_{1}>Z_{2}, \\
\eta =2  & \text{if } Z_{1}<Z_{2}. \\
  \end{cases}
$$
At the outset it is not clear that $\Delta^{*}$ is observable.
The following lemma shows that $\Delta^{*}$ is indeed a function
of $(\tilda{Z},\tilda{\Delta})$.

\begin{lem}
The following relation holds:
\[\{\Delta^{*}=1\}= \{\eta=1,\Delta_2=1\}\cup
\{\eta=2,\Delta_1=1\}\cup\{\eta=0,\Delta_1\vee\Delta_2=1\}.
\]
\end{lem}

\hop
Proof:
Suppose $\Delta^{*}=1$, i.e.~$T_1\wedge T_2\le C_1\wedge C_2$.
If $\eta=1$ then $T_1\wedge C_1> T_2\wedge C_2$.
We need to show that $\Delta_2=1$. Suppose the contrary,
i.e.~$T_2>C_2$; then $T_2\wedge C_2=C_2$ and we have
\[T_1\wedge C_1> C_2\arr T_1>C_2\]
which gives $T_1\wedge T_2 >C_2$ and hence
$T_1\wedge T_2> C_1\wedge C_2$, contradicting $\Delta^{*}=1$.
A similar argument with $\eta =2$ shows that the left hand side is
contained in the right hand side.

Next suppose $\eta=0,\Delta_{1}=1$,
i.e.~$T_1\wedge C_1=T_2\wedge C_2$ and $T_1\le C_1$.
These two easily yield $T_1\le C_1\wedge C_2$.
Now $\Delta_1=1$ implies $T_1\wedge C_1=T_1$
and together with $\eta=0$ gives
$T_1\wedge T_2=T_1$. Hence $\Delta^{*}=1$.
Similarly if $\eta=1,\Delta_2=1$,
i.e.~$T_1\wedge C_1>T_2\wedge C_2$
and $T_2\wedge C_2=T_2$. Then
$T_1\wedge C_1> T_2$ so that $T_2< C_1$.
All this results in
\[T_2\le C_1\wedge C_2
   \quad {\rm implying} \quad
   T_1\wedge T_2\le C_1\wedge C_2
   \quad {\rm implying} \quad\Delta^{*}=1. \]
Similar arguments handle the other cases. \qed

\smallskip
Turning to the distribution of
$(\tilda{Z},\tilda{\Delta})$,
$P^{\tilda{Z},\tilda{\Delta}}$
can be written as a product,
\[P^{(Z^{*},\Delta^{*})}
   P^{\eta}(\cdot\midd(Z^{*},\Delta^{*}))P^{Z_{\eta}}
  (\cdot\midd (Z^{*},\Delta^{*}),\eta)\]
where $Z_0 = Z^*$. We will show next that $P$ can be recovered from
\[P^{(Z^{*},\Delta^{*})},
   \quad P^{\eta}(\cdot\midd (Z^{*},\Delta^{*}=1)),
   \quad P^{Z_{\eta}}(\cdot\midd (Z^{*},\Delta^{*}=1),\eta). \]

(1)
The relation between $T^{*}$ and $(Z^{*},\Delta^{*})$ is
obvious. This is just a one-dimensional censoring model where
$T^{*}$ is censored by the independent variable $C^{*}$.
This gives
\[P\{T^{*}=t\midd T^{*} \geq t\}=
\frac{P\{Z^{*}=t,\Delta^{*}=1\}}{P\{Z^{*}\geq t\}}. \]

(2)
For $t$ with $P\{T^{*}=t\}>0$,
\begin{align*}
P\{\epsilon=1\midd T^{*}=t\}&=P\{\eta=1\midd Z^{*}=t,\Delta^{*}=1\}, \\
P\{\epsilon=2\midd T^{*}=t\}&=P\{\eta=2\midd Z^{*}=t,\Delta^{*}=1\}, \\
P\{\epsilon=0\midd T^{*}=t\}&=P\{\eta=0\midd Z^{*}=t,\Delta^{*}=1\}.
\end{align*}
This follows from
\begin{align*}
P\{\epsilon=1\midd T^{*}=t\}
&=\frac{P\{\epsilon=1, T^{*}=t\}}{P\{T^{*}=t\}} \\
&=\frac{P\{\epsilon=1, T_{1}>t, T_{2}=t, C_{1} \geq t, C_{2}\geq t\}}
{P\{T^{*}=t, C^{*}\geq t\}}\\
&=\frac{P\{\eta=1, Z^{*}=t,\Delta^{*}=1\}} {P\{Z^{*}=t,\Delta^{*}=1\}};
\end{align*}
the other cases are handled similarly.

(3)
For any $E$, let $P_{E}$ stand for the conditional
probability given $E$. For $t > j$,
\[P_{T^{*}=j,  \epsilon =1}\{T_{1}=t\midd T_{1}\geq t\}
  = \frac{P_{Z^{*}=j,\Delta^{*}=1,\eta =1}
\{Z_1=t, \Delta_1=1\}}{P_{Z^{*}=j,\Delta^{*}=1,\eta =1}
\{Z_1\geq t\}}. \]
To put it differently, what we have is the
following one-dimensional censoring model. Consider random
variables $\hat T$ distributed as $P_{T^{*}=j,\epsilon =1}$,
and $\hat C$ independent of $\hat T$ distributed as
$P^{C_2}\{\cdot\midd C_2>j\}$.
Then the pair $\hat Z=\hat T\wedge\hat C$,
$\hat\Delta=I\{\hat T\le \hat C\}$ has the same distribution
as $(Z_1,\Delta_1)$ under $P_{Z^{*}=j,\Delta^{*}=1}$.
The details are as follows:
\begin{align*}
P_{T^{*}=j,\epsilon=1}\{T_1=t\midd T_1\ge t\}
&=\frac{P\{T_1=t,T^{*}=j,\epsilon=1\}}
   {P\{T_1\ge t,T^{*}=j,\epsilon=1\}}\\
&=\frac{P\{T_1=t,T^{*}=j,C_1\geq t, C_2\ge j\}}
  {P\{T_1\ge t,T^{*}=j,C_1\geq t, C_2\ge j\}} \\
&=\frac{P\{Z_{1}=t,\Delta_1=1,Z^{*}=j,\Delta^{*}=1,\eta=1\}}
 {P\{Z_1\ge t,\Delta_1=1,Z^{*}=j,\Delta^{*}=1,\eta=1\}}\\
&=\frac{P_{Z^{*}=j,\Delta^{*}=1,\eta=1}
  \{Z_1=t, \Delta_1=1\}}{P_{Z^{*}=j,\Delta^{*}=1,\eta=1}
  \{(Z_1\ge t,\Delta_1=1\}}.
\end{align*}

(4)
Similarly,
\[P_{T^{*}=j,\epsilon=2}\{T_2=t\midd T_1\ge t\}
   = \frac{P_{Z^{*}=j,\Delta^{*}=1,\eta =2}
  \{Z_2=t,\Delta_2=1\}}{P_{Z^{*}=j,\Delta^{*}=1,\eta=2}
  \{Z_2\ge t,\Delta_2=1\}}. \]

We have thus obtained the distribution of $\tilda{T}$ from
that of $ (\tilda{Z}, \tilda{\Delta})$. A couple of
aspects of this identification deserves mention.

(a)
In order to recover the distribution of $T^*$, as
done in the Kaplan--Meier estimate, we use
the distribution of $(Z^*, \Delta^*)$.

(b)
In obtaining the conditional distributions
$P\{T_1 =t\midd T^*=t^*,\epsilon =1\}$
we only use factors of the form
$P\{\cdot\midd Z^*=t^*,\Delta^{*}=1\}$.
More specifically terms like
$P\{\cdot\midd Z^*=t^*,\Delta^{*}=0\}$ are not needed.
As mentioned earlier in this section, this amounts
to showing that the map
$\phi\colon F \mapsto F^{(\tilda{Z}, \tilda{\Delta})}$
is not onto. This phenomenon reappears
when we attempt to write expressions for the likelihood function.

The difficulty in the bivariate case stems from
fact that, in general, the empirical distributions is not in the
range of $\phi$ and consequently no unique inverse.
For instance let $Z_{1}^{*},\ldots,Z_{k}^{*}$ be observations
with $\Delta^{*}_i=1,\eta_i=1,\Delta_{1,i}=1$.
Suppose $Z_{0}$ is an observation with $\Delta^{*}=0,\eta_0=1$
and $\Delta_{0,1}=0,\Delta_{0,2}=1$. Assume that $Z_{0,1},Z_{0,2}$
are distinct from the rest. From the one-dimensional case
it follows that $P^{*}$ will give positive values to
$Z_{1}^{*},\ldots,Z_{k}^{*}$ and
$P\{T_1=Z_{0,1}\midd T_2=Z_i^{*}\}=0$
for $i=1,\ldots,k$. On the other hand,
since $\Delta^{*}=0$ corresponds to $T_2>Z_{1,0}^{*}$,
$P\{T_1=Z_{0,1}\midd T_2=Z_i^{*}\}>0$ for some $i$.

\subsection{The likelihood}

The observed random variable has three components, namely
$(Z^{*},\Delta^{*})$, $\eta$ and $(Z_{\eta},\Delta_\eta)$.
Consequently  the likelihood also has three factors
\[P^{(Z^{*},\Delta^{*})}, \quad
  P^{\eta}(\cdot\midd(Z^{*},\Delta^{*})), \quad
  P^{(Z_{\eta},\Delta_{\eta})}
  (\cdot\midd (Z^{*},\Delta^{*}),\eta). \]

(1)
The distribution of $(Z^{*},\Delta^{*})$ is the same as that
coming from the one-dimensional problem where the survival time is
$T^{*}$ and the  censoring variable is $C^{*}$.
The exact form of $P^{(Z^{*},\Delta^{*})}$ does not concern
us here. We note that $P^{(Z^{*},\Delta^{*})}$ depends only on $P^{*}$.

(2)
We have
\[\frac{P\{\eta=j,(Z^{*}=z,\Delta^{*}=1)\}}
{P\{Z^{*}=z,\Delta^{*}=1\}}
=\frac{P\{\eta=j,T^{*}=z,C^{*}\geq z\}}
{P\{Z^{*}=z,C^{*}\geq z\}}
=\frac{P\{\epsilon=j,T^{*}=z\}}
{P\{T^{*}=z\}}. \]
Hence
\[P\{\eta=j\midd Z^{*}=z,\Delta^{*}=1\}
   =P\{\epsilon=j\midd t^{*}=z\}. \]

(3)
We next move to the case
$P^{(Z_{\eta},\Delta_{\eta})}\{\cdot\midd Z^{*}, \Delta^{*}=1,\eta\}$.
This again corresponds to a one-dimensional model.
Let $\hat T$ have the distribution
$P^{T_{1}}\{\cdot\midd T^{*}=z^{*},\epsilon =1\}$
and $\hat C$ be independent of $\hat T$ have the distribution
$P^{C_{2}}\{\cdot\midd C_1>z^{*},C_2\ge z^{*}\}$.
Then given
$Z^{*}=z^{*},\Delta^{*}=1,\eta =1$, $(Z_{1},\Delta_{1})$
has the same distribution as $\hat T\wedge \hat C$,
$I\{\hat T\le \hat C\}$. More explicitly,
\begin{setlength}{\multlinegap}{0pt}
\begin{multline*}
\frac{P\{Z_{1}=z,\Delta_{1}=\delta_{1},Z^{*}=z^{*},\Delta^{*}=1,\eta=1\}}
   {P\{Z^{*}=z^{*}\}} \\
=\frac{P\{T_1=z, C_1\ge i, T_2=z^{*}, C_2\ge z^{*})\}}
  {P\{T_1>z^{*}, T_2=z^{*}, C_1>z^{*}, C_2\ge z^{*}\}} \\
=\frac{P\{T_1=z, T_2=z^{*}\}}{P\{T_1>z^{*}, T_2=z^{*}\}}
  \frac{P\{C_1\ge z, C_2\ge z^{*}\}}{P\{C_1>z^{*},C_2\ge z^{*}\}}.
\end{multline*}
\end{setlength}
The $\Delta_{1}=0$ cases are handled similarly.

What is left are the terms with $\Delta^{*}=0$, namely
$$P^{\eta}\{\cdot\midd\Delta^{*}=0\}, \quad
  P^{(Z_{\eta},\Delta_{\eta})}\{\cdot\midd\Delta^{*}=0,\eta\}. $$
The first of these depends only on $G$ and
the second does not have an explicit expression
in terms of the parametrization of $P$ that we have used.

\section{Prior and Posterior}

\noindent
We first look at the time-discrete case, where $P$ and $G$
are distribution with supports in
$\{1,2,3,\ldots\}\times\{1,2,3,\ldots\}$.
A prior $\Pi$ for $P$ is called a Beta process, see Hjort (1990), if
\begin{enumerate}[(a)]
\item the hazards $\{P\{T^{*}=i\midd T^{*}\geq i\}$
are independent and are distributed
as Beta random variables with parameters
$(\alpha^{*}(i),\beta^{*}(i))$;

\item $P^{\epsilon}\{\cdot\midd T^{*}=t\}$
has a Dirichlet distribution with parameters
$(\alpha(0\midd t),\alpha(1\midd t),\allowbreak \alpha(2\midd t))$;

\item for $\epsilon=1,2$, the hazards
$$P\{T_{\epsilon}=i+d\midd T^{*}=i,\epsilon\}/
  P\{T_{\epsilon}\ge i+d\midd T^{*}=i,\epsilon\} $$
are independent Betas with parameters
$(\alpha_{i,\epsilon}(d),\beta_{i,\epsilon}(d))$;

\item the families appearing in (a), (b), (c) are independent.
\end{enumerate}

With i.i.d.~observations $(\tilda Z_1,\tilda\Delta_1),\ldots,
(\tilda Z_n,\tilda\Delta_n)$, the likelihood factors as follows:
\begin{enumerate}[(a)]
\item Probabilities of $(Z_{1}^{*},\Delta_{1}^{*}),
\ldots,(Z_{n}^{*},\Delta_{n}^{*})$.
These arise from a one-dimensional censoring model involving
$P^{*}$ and $G$.
\item Probabilities of $\eta$ from $Z_{i}^{*}$s with
$\Delta_{i}^{*}=1$.
\item Probabilities of $Z_{j}$ from $Z_{i}^{*}$s with
$\Delta_{i}^{*}=1$ and $\eta =j$ for $j=1,2$.
\item probabilities arising from $Z_{i}^{*}$s with
$\Delta_{i}^{*}=0$. These do not have an explicit expression in
terms of the parameters $P^{*}$, $P^{\epsilon}\{\cdot\midd T^{*}\}$, etc.
\end{enumerate}
In order to arrive at a tractable `posterior' we will ignore
terms in `(d)' above and work with an {\it incomplete likelihood}
involving only (a), (b), (c). We view the terms related 
to (a), (b), (c) as more statistically relevant, 
for the estimation of the bivariate survival curve, 
than those related to (d). If this is done then it is easy
to see that the resulting `posterior' is again a Beta process.
Somewhat more formally, assume for simplicity that
$\eta_{i}\neq 0$ for $i=1,\ldots,n$, and set
\begin{align*}
Y_{i}^{*} &=\#\{j\colon Z_j^{*}\ge i\}, \\
\Delta N_i^{*} &= \# \{j\colon Z_j^{*} = i, \Delta^*_j=1\}, \\
N_i(1) &=\#\{j\colon Z_j^{*} = i, \Delta^{*}_j=1,\eta=1\}, \\
N_i(2) &=\#\{j\colon Z_j^{*} = i, \Delta^{*}_j=1,\eta=2\}, \\
  Y_{i,d}(1) &= \# \{j\colon Z_j^{*} = i,
  \Delta^{*}_j=1,\eta=1,Z_{1,j}\geq i+d\}, \\
\Delta N_{i,d}(1) &= \# \{j\colon Z_j^{*} = i,
  \Delta^{*}_j=1,\eta=1,Z_{1,j}= i+d,\Delta_{1,j}=1\}, \\
Y_{i,d}(2) &= \# \{j\colon Z_j^{*} = i,
  \Delta^{*}_j=1,\eta=2,Z_{2,j}\geq i+d\}, \\
\Delta N_{i,d}(1) &= \# \{j\colon Z_j^{*} = i,
\Delta^{*}_j=1,\eta=2,Z_{2,j}= i+d,\Delta_{2,j}=1\}.
\end{align*}
The updated parameters are
\begin{align*}
\hat{\alpha}^{*}(i)&=\alpha^{*}(i)+\Delta N_i^{*}, \\
\hat{\beta}^{*}(i)&=\beta^{*}(i))+ Y_{i}^{*}-\Delta N_i^{*}, \\
(\hat\alpha(0\midd i),\hat\alpha(1\midd i),\hat\alpha(2\midd i))
&=(\alpha(0\midd t)+N_i(0),\alpha(1\midd t)+N_i(1),
   \alpha(2\midd t)+N_i(2)), \\
\hat\alpha_{i,1}(d)&=\alpha_{i,1}(d)+\Delta N_{i,d}(1), \\
\hat\beta_{i,1}(d)&=\beta_{i,1}(d))+Y_{i,d}(1)-\Delta N_{i,d}(1), \\
\hat\alpha_{i,2}(d)&=\alpha_{i,2}(d)+ \Delta N_{i,d}(2), \\
\hat\beta_{i,2}(d)&=\beta_{i,2}(d)+Y_{i,d}(2)-\Delta N_{i,d}(2). \\
\end{align*}
An easy computation gives the Bayes estimates as
\begin{align*}
\E_n P\{T^{*}=i\midd T^{*}\geq i\}
  &=\frac{\alpha^{*}(i)+\Delta N_i^{*}}
  {\alpha^{*}(i)+\beta^{*}(i)+Y_{i}^{*}}, \\
  \E_n P\{T^{*}>i\midd T^{*}\geq i\}
  &=1-\frac{\alpha^{*}(i)+\Delta N_i^{*}}
  {\alpha^{*}(i)+\beta^{*}(i)+Y_{i}^{*}} \\
  &= \frac{\beta^{*}(i)+Y_{i}^{*}-\Delta N_i^{*}}
  {\alpha^{*}(i)+\beta^{*}(i)+Y_{i}^{*}}, \\
\E_n P^{\epsilon}\{j\midd T^{*}=i\}
   &=\frac{\alpha(j\midd t)+N_i(j)}
   {\alpha(0\midd t)+N_i(0)+
   \alpha(1\midd t)+N_{i}(2)+\alpha(2\midd t)+N_{i}(2)}, 
\end{align*}
\begin{align*}
\E_n\frac{P\{T_{1}=i+d\midd T^{*}=i,\epsilon=1\}}
   {P\{T_{1}\ge i+d\midd T^{*}=i,\epsilon=1\}}
   &=\frac{\alpha_{i,1}(d)+\Delta N_{i,d}(1)}
   {\alpha_{i,1}(d)+\beta_{i,1}(d)+Y_{i,d}(1)},\\
   \E_n\frac{P\{T_{1}>i+d\midd T^{*}=i,\epsilon=1\}}
   {P\{T_{1}\ge i+d\midd T^{*}=i,\epsilon=1\}}
   &=1-\frac{\alpha_{i,1}(d)+\Delta N_{i,d}(1)}
   {\alpha_{i,1}(d)+\beta_{i,1}(d)+Y_{i,d}(1)} \\
  &=\frac{\beta_{i,1}(d)+Y_{i,d}(1)-\Delta N_{i,d}(1)}
   {\alpha_{i,1}(d)+\beta_{i,1}(d)+Y_{i,d}(1)}.
\end{align*}
We also have
\begin{setlength}{\multlinegap}{0pt}
\begin{multline}\label{eqa}
P\{T_{1}=i+d, T_{2}=i\}
=\Bigl[\prod_{1}^{i-1}\frac{P\{T^{*}>j\}}{P\{T^{*}\ge j \}}\Bigr]
\frac{P\{T^{*}=i\}}{P\{T^{*}\ge i\}}
   P\{\epsilon=1\midd T^{*}=i\} \\
\Bigl[\prod_{1}^{d-1}\frac{P\{T_1 >i+d'\midd T^{*}=i,
   \epsilon=1\}} {P\{T_1\ge i+d'\midd T^{*}=i,\epsilon =1\}}\Bigr]
\Bigl[\frac{P\{T_1=i+d\midd T^{*}=i,\epsilon =1\}}
{P\{T_1\ge i+ d\midd T^{*}=i,\epsilon=1\}}\Bigr].
\end{multline}
\end{setlength}

As we let the prior parameters go to zero, for each
$j$, the Bayes estimate of 
$P\{T^{*}>j\}/P\{T^{*}\ge j\}$ converges to 
$N_{j}^{*+}/N_j^{*}$ where 
$N_{j}^{*+}=Y_{j}^{*}-\Delta N_{j}^{*}$ and the Bayes estimate of
$P\{T^{*}=i\}/P\{T^{*}\ge i\}$ converges to  
$\Delta N_{i}^{*}/Y_{i}^{*}$. Carrying out a similar evaluation 
of the limit in each of the factors in (\ref{eqa}) and setting
$N_{i,d'}^{+}(1)= Y_{i,d'}(1)-\Delta N_{i,d'}(1)$ we get the
`noninformative' Bayes estimator
\[ \Bigl\{\prod_{1}^{i-1}\frac{N_{j}^{*+}}{N_j^{*}}\Bigr\}
\frac{\Delta N_{i}^{*}}{Y_{i}^{*}} \frac{N_{i}(1)}{\Delta N_{i}^{*}}
\Bigl\{\prod_{1}^{d-1}
\frac{N_{i,d'}^{+}(1)}{N_{i,d'}(1)}\Bigr\}
\frac{\Delta N_{i,d}(1)}{ Y_{i,d}(1)}. \]
While we have not worked out the details, 
we do expect this noninformative Bayes estimator 
to be consistent.

\section{Beta processes}

\noindent
In the context of general Beta processes it is convenient to work
with the cumulative hazard function of a distribution rather than the
distribution function itself. These processes were introduced in
Hjort (1990). Recall that if $F$ is a distribution function and if
$\bar F=1-F$ is the survival function then the cumulative hazard
function $A$ is the function
\[A(t)=\int_{[0,t]}\frac{\d F(s)}{F[s,\infty)}. \]
The $A$ is right continuous and increasing,
has the same jump points as $F$, and if $t$
is a jump point then the jump size
is $\Delta A(t)=\Delta F(t)/F[t,\infty)$.

Let $F_{0}$ be a distribution function and let $A_{0}^d,A_{0}^{c}$
be the discrete and continuous parts of its cumulative
hazard function $A_{0}$. Let $w$ be a continuous positive function
on $\mathbb{R}^{+}$. A Beta process with parameters $w$ and
$A_0$ is a process such that the process $\{A(t)\colon t>0\}$
has independent increments and has the L\'evy measure
\[\d L(t,s)=w(t)s^{-1}(1-s)^{w(t)-1}\,\d A^{c}_0(t)\,\d s
   \quad {\rm for\ }t\in\mathbb{R}^{+},\,0\leq s \leq 1, \]
and if $t_0$ is a jump point then
\[\Delta A(t_0)\sim\text{Beta}\{w(t_0)A^{d}_0(t_0),
   w(t_0)(1-A^{d}_{0}(t_0))\}. \]
If we have $n$ i.i.d.~observations with some of them right censored
and if we set
\[\Delta N(t)= \#\{i\colon Z_i=t,\Delta_i=1\}, \quad
   Y(t)= \#\{i\colon Z_i \ge t\}, \]
then Hjort (1990) has shown that the posterior
is again a Beta process with parameters
$w+Y$ and $\int_{0}^{\cdot}(w\,\d A_0+\Delta N)/(w+Y)$.

Since our approach to the bivariate problem is a sequence of
one-dimensional censoring models, a natural generalization
of the one-dimensional Beta process suggests itself.
Let $F_0$ be a distribution function on
$\mathbb{R}^{+}\times\mathbb{R}^{+}$.
Let $F_{0}^{*}$ be the distribution of $T^{*}$,
$F^{\epsilon}_{0}(\cdot\midd T^{*})$ be the distribution of $\epsilon$
given $T^{*}$ and $F_{0}^{1}(\cdot\midd T^{*})$ the distribution of
$T_{1}$ given $(T^{*}, \epsilon =1 )$ and $F_{0}^{2}(\cdot\midd T^{*})$
the distribution of $T_{2}$ given $\{T^{*}, \epsilon =2)\}$.
Let next $w$ be a function on $\mathbb{R}^{+}\times \mathbb{R}^{+}$.
We write $w^{*}(t)=w(t,t)$. Furthermore if $s>t$, then
set $w_{1}(s\midd t) = w(t,s)$. Define  $w_{2}(s\midd t)$ similarly.

Let $\alpha_{0}$ be a function from $\mathbb{R}$ to
$(\mathbb{R})^{3}$. The parameter space is

(1)
$\mathcal{M}$ -- the set of all distributions on
$\mathbb{R}^{+}$. We view this as the set of all possible
distributions for $T^{*}$.

(2)
$\Gamma$ -- all measurable functions from $\mathbb{R}$
to the set of probability measures on $\{0,1,2\}$.
This correspond to the set of conditional distributions
of $\epsilon$ given $T^{*}$.

(3)
$\Gamma_1$ -- all measurable functions from $\mathbb{R}$
to $\mathcal{M}$, viewed as the set of distributions
for $T_1$ given $T^{*}$ and $\epsilon=1$.

(4)
$\Gamma_{2}$ -- all measurable functions from $\mathbb{R}$
to $\mathcal{M}$, viewed as the set of distributions for
$T_{2}$ given $T^{*}$ and $\epsilon=2$.

A prior $\Pi=\Pi_1\times\Pi_2\times\Pi_3\times\Pi_4$
on $\mathcal{M}\times\Gamma\times\Gamma_1\times\Gamma_2$
is called a bivariate Beta process prior if

(1)
$\Pi_1$ is a Beta process with parameters $C^{*},A_{0}^{*}$
on $\mathcal{M}$;

(2) under $\Pi_2$, $\{\gamma(t^{*},\cdot)\colon t^{*} \in \mathbb{R}\}$
are independent and for each $t^{*}$, $\gamma(t^{*},\cdot)$
has a Dirichlet distribution with parameters
$(\alpha_0(t^{*}),\alpha_1(t^{*}),\alpha_2(t^{*})$;

(3)
under $\Pi_3$, $\{\gamma_{1}(t^{*},\cdot)\colon t^{*} \in \mathbb{R}\}$
are independent and for each $t^{*}$,
$\gamma_1(t^{*},\cdot)$ is a one-dimensional Beta process with
parameters $w_{1}(\cdot\midd t^{*})$ and $A_0^{(1)}(\cdot\midd t^{*})$,

(4)
under $\Pi_4$, $\{\gamma_{2}(t^{*},\cdot)\colon t^{*}\in\mathbb{R}\}$
are independent and for each $t^{*}$, $\gamma_{2}(t^{*},\cdot)$
is a one-dimensional Beta process with parameters
$w_{2}(\cdot\midd t^{*})$ and $A_0^{(2)}(\cdot\midd t^{*})$.

A technical argument is needed to properly demonstrate 
the existence of the above process; see the appendix.

It is fairly routine to see that given $n$ i.i.d.~observations,
using the incomplete likelihood, the posterior is again
a bivariate Beta process. The independence assumption
in the prior and the product form of the likelihood
ensures the required independence for the posterior.
That the posterior is again a Beta process follows from
the one-dimensional results. The exact Bayes estimators
can be written down but become notationally tedious.
We just note that if we let $w$ go to zero
we get exactly the same noninformative-prior estimator as before.


\section{An example}

\noindent
We first examine Pruitt's (1988) example via our approach.
The only observations with $\Delta^{*}=1$ are the ones with both
coordinates uncensored. Thus the estimate of the distribution of
$T^{*}$ will have positive mass only at $Z^{*}$ corresponding to
these points. The censored observations will all have $Z^{*}=0$.
Also it is only the observations with both coordinates uncensored
that are used to estimate the conditional probabilities.

Considering a simple numerical example in Dabrowska (1988),
suppose that the observed points written as
$(Z_{1},Z_{2},\Delta_{1},\Delta_{2})$ are
\[(.51,.02,1,1),\,(.11,.62,1,0),\,(.24,0.24,0,0),\,(.68,.68,1,1). \]
The noninformative-prior estimate will give mass 1/4 to each of
(.11,.62), (.51,.02) and 1/2 to the point (.68,.68).

\begin{tabular}{|c|c|c|}
  Estimate of  & Dabrowska & Noninformative-prior \\
\mbox{$P\{T_{1}>0,T_{2}>0\}$} & 1 & 1 \\
\mbox{$P\{T_{1}>.11,T_{2}>0\}$} & .75 & .75 \\
\mbox{$P\{T_{1}>.24,T_{2}>0\}$}  & .75 & .75 \\
\mbox{$P\{T_{1}>.51,T_{2}>0\}$}  & .375 & .5 \\
\mbox{$P\{T_{1}>.68,T_{2}>0\}$}  & 0 & 0 \\
\mbox{$P\{T_{1}>0,T_{2}>.02\}$} & .75 & .75 \\
\mbox{$P\{T_{1}>0,T_{2}>.24\}$} & .75 & .75 \\
\mbox{$P\{T_{1}>0,T_{2}>.62\}$} & .75 & .5 \\
\mbox{$P\{T_{1}>0,T_{2}>.68\}$} & 0& 0 \\
\mbox{$P\{T_{1}>x,T_{2}>y\}$}  & .5 & .5 \\ \hline
\end{tabular}

\medskip\noindent
In the last row $x\in [.11,.68]$ and $y\in[.02,.68]$.

The set in the first column of the fourth row contains 
that in the last row, for $x>.51$. However, using 
the Dabrowska estimator method, the set in the last row 
has larger probability, leading to negative mass for some sets. 
Pruitt (1991) has examined the occurrence of negative mass 
in the Dabrowska estimate. This lack of monotonocity 
will not occur with the noninformative estimator.

\section{Appendix}

\noindent
Consider $\Gamma_{1}$. This is a subset of
$\mathcal{M}^\mathbb{R}$. Kolmogorov's consistency theorem
guarantees the existence of a product measure on
$\mathcal{M}^\mathbb{R}$. However we need to show that the outer
measure of $\Gamma_{1}$ is 1. This follows from the following
lemma. A similar argument appears in Judd (1985).

\begin{lem} Let $\mathcal{X}$ and $\mathcal{Y}$ be a complete
separable metric spaces with their corresponding Borel
$\sigma$-fields. Let $\mu=\prod_{i\in \mathcal{Y}}\mu_{i}$ be a
product of probability measures on $\mathcal{X}^{\mathcal{Y}}$. If
$\Gamma \subset \mathcal{X}^{\mathcal{Y}}$ is the set of all
measurable functions from $\mathcal{Y}$ to $\mathcal{X}$, then the
$\mu$ outer measure of $\Gamma$ is one.
\end{lem}

\hop
Proof: 
Let $E$ be any measurable set in $M(\mathcal{Y})^\mathcal{X}$ 
containing $\Gamma$. We will show that $E$ is the whole space. 
Since $E$ is measurable there is a countable set
$\{t_1,t_2,\ldots\}$ such that $E$ is in the $\sigma$-algebra 
generated by the projections to these coordinates. 
Hence for any $\gamma_1,\gamma_2$, 
if $\gamma_1(\cdot |t_i)= \gamma_2(\cdot |t_i)$ for $i=1,2,\ldots$,
then either both $\gamma_1$ and $\gamma_2$ or none of them 
belong to $E$. Take any $\gamma\in M(\mathcal{Y})^\mathcal{X}$, 
and pick a $\gamma_0$ from $\Delta$. Define $\gamma^*$ by 
$$\gamma^*(\cdot\midd t_i)=
\begin{cases}
\gamma(\cdot\midd t_i) &\mbox{for\ } i=1,2,\ldots; \\
\gamma_0(\cdot\midd t) &\mbox{otherwise}. 
\end{cases}$$
Since $\gamma_0$ is modified only at countably many points, 
measurability still holds and hence $\gamma^* \in \Gamma$.
Now $\gamma^*(\cdot\midd t_i)= \gamma(\cdot\midd t_i)$ 
for $i=1,2,\ldots$, hence $\gamma$ is also in $E$. \qed 


\def\annals{Annals of Statistics}
\def\jasa{Journal of the American Statistical Association}


\bigskip
\centerline{\bf References}

\def\ref#1{{\noindent\hangafter=1\hangindent=20pt
  #1\smallskip}}
\parindent0pt
\baselineskip16pt
\parskip3pt
\medskip

\ref{%
[1] Andersen, P.K., Borgan, \O., Gill, R.D.~and Keiding, N. (1993).
{\sl Statistical Models Based on Counting Processes.}
Springer-Verlag, New York.}

\ref{%
[2] Dabrowska, D.M. (1988).
Kaplan--Meier estimate on the plane.
{\sl\annals} {\bf 15}, 1475--1489.}

\ref{%
[3] Hjort, N.L. (1990).
Nonparametric Bayes estimators based on Beta
processes in models for life history data.
{\sl\annals} {\bf 18}, 1259--1294.}

\ref{%
[4] Judd, K.L. (1985).
The law of large numbers with a continuum of
IID random variables.
{\sl Journal of Econometric Theory} {\bf 35}, 19--25.}

\ref{%
[5] Langberg, N.A.~and Shaked, M. (1982).
On the identifiability of multivariate life distribution functions.
{\sl Annals of Probability} {\bf 10}, 773--779.}

\ref{%
[6] Peterson, A.V. (1977).
Expressind the Kaplan--Meier estimator as a
function of empirical subsurvival functions.
{\sl\jasa} {\bf 72}, 854--858.}

\ref{%
[7] Pruitt, R.C. (1988).
An inconsistent Bayes estimate in bivariate
survival curve analysis.
Preprint 529, University of Minnestota.}

\ref{%
[8] Pruitt, R.C. (1991).
On negative mass assigned in the bivariate
Kaplan--Meier estimator.
{\sl\annals} {\bf 19}, 443--453.}

\ref{%
[9] Pruitt, R.C. (1993).
Identifiability of bivariate survival curves
from censored data.
{\sl\jasa} {\bf 88}, 573--579.}


\end{document}